\documentclass[12pt]{amsart}

\usepackage[all]{xy}
\usepackage{amscd}
\usepackage{amssymb}

\newtheorem{thm}{Theorem}[section]
\newtheorem{lem}[thm]{Lemma}

\newtheorem{prop}[thm]{Proposition}
\newtheorem{conj}[thm]{Conjecture}

\theoremstyle{definition}

\newtheorem{rem}[thm]{Remark}

\theoremstyle{remark}

\numberwithin{equation}{section}

\def\F{{\mathbb F}}

\def\Q{{\mathbb Q}}
\def\Z{{\mathbb Z}}
\def\C{{\mathbb C}}
\def\O{{\mathcal O}}

\def\X{{\mathfrak X}}
\def\Y{{\mathfrak Y}}
\def\ZZ{{\mathfrak Z}}

\def\Gr{\text{\rm Gr}}

\def\Gal{\text{\rm Gal}}

\def\GL{\text{\rm GL}}

\def\Spec{\text{\rm Spec}\,}
\def\sep{\text{\rm sep}}

\begin{document}

\title[Weight-monodromy conjecture]
{Weight-monodromy conjecture over equal
characteristic local fields}

\author[Tetsushi Ito]{Tetsushi Ito}
\address{Department of Mathematics, Faculty of Science,
Kyoto University, Kyoto 606-8502, Japan}
\email{tetsushi\char`\@math.kyoto-u.ac.jp}

\date{\today}

\begin{abstract}
The aim of this paper is to study certain properties of
the weight spectral sequences of Rapoport-Zink by a specialization argument.
By reducing to the case over finite fields
previously treated by Deligne, we prove that
the weight filtration and the monodromy filtration
defined on the $l$-adic \'etale cohomology coincide,
up to shift, for proper smooth varieties over equal characteristic local fields.
We also prove that the weight spectral sequences degenerate at $E_2$
in any characteristic without using log geometry.
Moreover, as an application, we give a modulo $p>0$ reduction proof
of a Hodge analogue previously considered by Steenbrink.
\end{abstract}

\maketitle

\section{Introduction}
\label{SectionIntroduction}

The aim of this paper is to study certain properties of
the weight spectral sequences of Rapoport-Zink by a specialization argument.
Let $K$ be a discrete valuation field with ring of integers
$\O_K$ and residue field $F$. Assume that $\O_K$ is henselian.
Let $\pi$ be a uniformizer of $K$,
and $l$ a prime number different from the characteristic of $F$.
Let $\X$ be a regular scheme proper and flat over $\O_K$
such that the generic fiber $\X_K := \X \otimes_{\O_K}\!K$
is a proper smooth variety of dimension $n$ over $K$,
and the special fiber $\X_F := \X \otimes_{\O_K}\!F$ is
a divisor of $\X$ with simple normal crossings.
Such $\X$ is called a {\em proper strictly semistable scheme} over $\O_K$.
Let $\X^{(k)}$ be the disjoint union of $k$ by $k$ intersections of
the irreducible components of $\X_F$.

The {\em weight spectral sequence} of Rapoport-Zink
is the following $\Gal(K^{\sep}/K)$-equivariant spectral sequence
relating the $l$-adic \'etale cohomology of
the special fiber and the generic fiber
(\cite{Rapoport-Zink}, Satz 2.10):
$$
   E_1^{-r,\,w+r} = \bigoplus_{k \geq \max\{0,-r\}}
   H^{w-r-2k}_{\text{\rm \'et}}
     \big( \X^{(2k+r+1)}_{F^{\sep}},\, \Q_l(-r-k) \big)
   \ \Rightarrow \ H^w_{\text{\rm \'et}}(\X_{K^{\sep}},\Q_l),
$$
where $\X^{(k)}_{F^{\sep}} := \X^{(k)} \otimes_F\!F^{\sep}$,
\ $\X_{K^{\sep}} := \X_{K} \otimes_K\!K^{\sep}$,
and $F^{\sep}$ (resp.\ $K^{\sep}$) denotes
the separable closure of $F$ (resp.\ $K$).
We have a natural map, called the {\em monodromy operator}
$N \colon E_1^{i,\,j}(1) \to E_1^{i+2,\,j-2}$ such that
$N^r \colon E_1^{-r,\,w+r}(r) \xrightarrow{\cong} E_1^{r,\,w-r}$
is an isomorphism for all $r,w$.
Note that, the monodromy operator $N$ is a combination of
identity maps and zero maps,
and the differentials $d_1^{i,\,j}$ on the $E_1$-terms are
explicitly described
in terms of restriction morphisms and Gysin morphisms
(for details, see \cite{Rapoport-Zink}, \cite{Illusie1},
\cite{Illusie2}, \cite{TSaito}).

The main theorem of this paper is as follows.

\begin{thm}
\label{MainTheorem}
Let notation be as above.
\begin{enumerate}
\item {\rm (\cite{CNakayama}, Theorem 0.1)}
      The weight spectral sequence degenerates at $E_2$.
\item If $K$ is of equal characteristic, the monodromy operator $N$
induces an isomorphism
$N^r \colon E_2^{-r,\,w+r}(r) \xrightarrow{\cong} E_2^{r,\,w-r}$
for all $r,w$.
\end{enumerate}
\end{thm}

When $F$ is a finite field, the first part of Theorem \ref{MainTheorem}
is an easy consequence of the Weil conjectures
(\cite{Rapoport-Zink}, Satz 2.10, see also 
Proposition \ref{Degeneracy_weak}).
The general case was already proved by C.\ Nakayama
by specializing log schemes (\cite{CNakayama}, Theorem 0.1).
He essentially used his construction of
the weight spectral sequences for log schemes.
In this paper, we give another proof
without using log geometry.

The second part of Theorem \ref{MainTheorem}
is sometimes called the {\it weight-monodromy conjecture}.
Namely, we have the following conjecture in any characteristic
(\cite{HodgeI}, \cite{WeilII},
see also \cite{Rapoport-Zink}, \cite{Illusie1}, \cite{Illusie2}).

\begin{conj}[Weight-monodromy conjecture]
\label{WMC}
Let notation be as above. Then, the monodromy operator $N$
induces an isomorphism
$N^r \colon E_2^{-r,\,w+r}(r) \xrightarrow{\cong} E_2^{r,\,w-r}$
for all $r,w$.
\end{conj}

According to Theorem \ref{MainTheorem},
Conjecture \ref{WMC} holds if $K$ is of equal characteristic.
In mixed characteristic,
Conjecture \ref{WMC} is known to hold
if either ${\X}_K$ is an abelian variety or $\dim {\X}_K \leq 2$
(\cite{SGA7-I}, \cite{Rapoport-Zink}, \cite{deJong}, \cite{TSaito}).
Although some partial results were obtained (\cite{Ito2}, \cite{Ito3}),
Conjecture \ref{WMC} is still open in general
(see also Remark \ref{KnownCasesWMC}).

It is well-known that, if the residue field $F$ is finite,
Conjecture \ref{WMC} is equivalent to another conjecture
(Conjecture \ref{WMC_filtration})
on the monodromy filtration and the weight filtration
defined on the $l$-adic \'etale cohomology
of the geometric generic fiber $\X_{K^{\sep}}$
(Proposition \ref{WMC_equivalence}, Remark \ref{RelationTwoConjectures}).
When $\X/\O_K$ is the henselization of a family of varieties
over a curve over a finite field,
Deligne proved Conjecture \ref{WMC_filtration}, and hence
he essentially proved Conjecture \ref{WMC}
although the weight spectral sequence
was not constructed at that time
(Remark \ref{KnownCasesWMC}, \cite{WeilII}, Th\'eor\`eme 1.8.4).

In this paper, we give a proof of Theorem \ref{MainTheorem}
by a specialization argument.
The key geometric construction is given in \S \ref{SectionConstruction}
using N\'eron's blowing up.
Then, the first part of Theorem \ref{MainTheorem} easily follows
in all characteristics (Proposition \ref{PropDegeneracy}).
For the second part of Theorem \ref{MainTheorem},
we reduce the general case to the case treated by Deligne in \cite{WeilII}.
This argument is valid only in the equal characteristic case
(Remark \ref{MixedCharacteristicCase}).
Note that the possibility of such an argument was pointed out by
Illusie in \cite{Illusie2}, 8.7,
and a similar specialization argument for
Conjecture \ref{WMC_filtration} was considered by
Terasoma when $F$ is a finite field (\cite{Terasoma}).
In \S \ref{SectionConsequence}, we also give some consequences.
In \S \ref{SubsectionApplicationEqualChar}, we prove
the equal characteristic case of Conjecture \ref{WMC_filtration}
by a similar specialization argument
(Proposition \ref{EqualCharWMC}).
Moreover, in \S \ref{SubsectionApplicationHodge},
we give a modulo $p>0$ reduction proof
of a Hodge analogue previously considered by Steenbrink
(Proposition \ref{HodgeAnalogue}).

\vspace{0.2in}

\noindent
{\bf Acknowledgments.}
This paper is a revised version of the Master's thesis
of the author at University of Tokyo on March 2001 (\cite{Ito1}).
He is grateful to his advisor Takeshi Saito for his advice
and support during the graduate course.
He would like to thank Kazuya Kato for his support and
encouragement, Luc Illusie for invaluable discussion and comments
on the weight spectral sequences,
and Kazuhiro Fujiwara for the suggestion to use
N\'eron's blowing up.
He would also like to thank Yoichi Mieda for
reading a manuscript and comments.
The author was supported by the Japan Society for the
Promotion of Science Research Fellowships for Young
Scientists.

\section{Monodromy filtration and weight filtration}
\label{SectionFiltrations}

In this section, we recall the monodromy filtration and
the weight filtration defined
on the $l$-adic \'etale cohomology of proper smooth varieties
over local fields.
We also recall a conjecture on these two filtrations
(Conjecture \ref{WMC_filtration})
which is equivalent to Conjecture \ref{WMC} in the semistable reduction case.

As in \S \ref{SectionIntroduction},
let $K$ be a discrete valuation field with ring of integers
$\O_K$ and residue field $F$. Assume that $\O_K$ is henselian.
Let $\pi$ be a uniformizer of $K$,
and $l$ a prime number different from the characteristic of $F$.
Let $X_K$ be a proper smooth variety of dimension $n$ over $K$,
$V := H^w_{\text{\rm \'et}}(X_{K^{\sep}},\Q_l)$
the $l$-adic \'etale cohomology of
$X_{K^{\sep}} := X_K \otimes_K K^{\sep}$,
and $\rho \colon \Gal(K^{\sep}/K) \to \GL(V)$ the action
of $\Gal(K^{\sep}/K)$ on $V$.

\subsection{Monodromy filtration}
\label{SubsectionMonodromyFiltration}

First of all, recall that we have the exact sequence
$1 \to I_K \to \Gal(K^{\sep}/K) \to \Gal(F^{\sep}/F) \to 1$,
where $I_K$ is called the {\it inertia group} of $K$.
The pro-$l$-part of $I_K$ is isomorphic to $\Z_l(1)$ by a natural map
$$ t_l \colon I_K \ni \sigma \mapsto
   \left( \frac{\sigma(\pi^{1/l^m})}{\pi^{1/l^m}} \right)_{\!\!m}
   \in \varprojlim \mu_{l^m} =: \Z_l(1), $$
where $\mu_{l^m}$ is the group of $l^m$-th roots of unity
in $K^{\sep}$ (\cite{Serre}).
By Grothendieck's monodromy theorem (\cite{SerreTate}, Appendix,
see also \cite{SGA7-I},\ I, Variante 1.3),
there is an open subgroup $J \subset I_K$ and a nilpotent map
$N \colon V(1) \to V$ such that
  $\rho(\sigma) = \exp(t_l(\sigma) N)$
for all $\sigma \in J \subset I_K$.
Here $N$ is a nilpotent map means that
there is $r \geq 1$ such that $N^r \colon V(r) \to V$ is a zero map,
where $V(r)$ denotes the $r$-th Tate twist of $V$.
The map $N$ is called the {\it monodromy operator} on $V$.

The {\it monodromy filtration} $M$ on $V$
is a unique increasing filtration
such that $M_{-k} V = 0,\ M_k V = V$ for sufficiently large $k$,
$N (M_k V(1)) \subset M_{k-2} V$ for all $k$,
and $N$ induces an isomorphism
$N^k \colon \Gr^M_{k} V (k) \xrightarrow{\cong} \Gr^M_{-k} V$
for all $k \geq 0$, where $\Gr^M_{k} V := M_{k} V / M_{k-1} V$
(\cite{WeilII}, I, (1.7.2)).

\subsection{Weight filtration}
\label{SubsectionWeightFiltration}

Usually, the notion of weights is considered when
$F$ is a finite field.
Here we consider it in a slightly general situation.
Assume that $F$ is a purely inseparable extension of
a finitely generated extension of a prime field (i.e.\ $\F_p$ or $\Q$).
It is always possible to find a finitely generated $\Z$-algebra $A$
contained in $F$ such that $F$ is a purely inseparable extension
of the field of fractions $\text{\rm Frac}\,A$ of $A$.

We say a continuous $l$-adic representation $\rho$ of
$\Gal(F^{\sep}/F)$ has {\it weight} $k$
if $\rho$ comes from a smooth $l$-adic sheaf $\mathcal F$ on
an open dense set $U \subset \Spec A$ by base change
and $\mathcal F$ has weight $k$ in the usual sense
(for example, see \cite{HodgeI}, \cite{WeilII}, (1.2)).
Namely, for all closed points $s \in U$,
the eigenvalues of the geometric Frobenius element at $s$
acting on ${\mathcal F}_{\bar{s}}$ are algebraic integers
such that the complex absolute values of 
their complex conjugates are $|\kappa(s)|^{k/2}$,
where $|\kappa(s)|$ denotes the order of the residue field
$\kappa(s)$ at $s$.
It is easy to see that this definition is independent of
the choice of $A,\ U,\ {\mathcal F}$.
By the Weil conjectures (\cite{WeilI}, \cite{WeilII}),
for a proper smooth variety $Z$ over $F$,
the $l$-adic \'etale cohomology $H^k_{\text{\rm \'et}}(Z_{F^{\sep}},\Q_l)$
has weight $k$ in this sense.

Let us go back to the situation in the beginning of this section.
As above, assume that $F$ is a purely inseparable extension of
a finitely generated extension of a prime field.
The {\it weight filtration} $W$ on $V$ 
is a unique increasing filtration such that
$W_{-k} V = 0,\ W_k V = V$ for sufficiently large $k$,
the action of $I_K$ on $\Gr^W_{k} V := W_{k} V / W_{k-1} V$
factors through a finite quotient,
and, after replacing $K$ by a finite extension of it,
$\Gal(F^{\sep}/F)$ acts on $\Gr^W_{k} V$ and
this action has weight $k$ (\cite{HodgeI}, \cite{WeilII}, I, (1.7.5),
for the existence of the weight filtration,
see \S \ref{SubsectionApplicationWeightSpectralSequence} below).

\subsection{An application to the weight spectral sequences}
\label{SubsectionApplicationWeightSpectralSequence}

Let notation be as in the beginning of this section.
Assume that the residue field  $F$ is
a purely inseparable extension of
a finitely generated extension of a prime field,
and there exists a proper strictly semistable
scheme $\X$ over $\O_K$ such that the generic fiber
$\X_K$ is isomorphic to $X_K$.
Then, each $E_1^{-r,\,w+r}$ of the weight spectral sequence
has weight $(w-r-2k)-2(-r-k) = w+r$ by the Weil conjectures
(\S \ref{SubsectionWeightFiltration}).
Therefore, the filtration on
$V := H^w_{\text{\rm \'et}}(X_{K^{\sep}},\Q_l)$
induced by the weight spectral sequence is nothing but
the weight filtration in \S \ref{SubsectionWeightFiltration}.
This proves the existence of the weight filtration on $V$
in the semistable reduction case.
In general, we may use de Jong's alteration
to reduce to the semistable reduction case
(Remark \ref{RemarkEqualChar}, see also the Introduction of \cite{deJong}.
See \cite{GoertzHaines} for the construction of the weight spectral sequence
via the weight filtration on nearby cycles).

Here we give an immediate application of the notion of weights
to the weight spectral sequences.

\begin{prop}
\label{Degeneracy_weak}
Let notation be as in \S \ref{SectionIntroduction}.
Assume that the residue field $F$ is a purely inseparable extension of
a finitely generated extension of a prime field.
Then, the weight spectral sequence degenerates at $E_2$.
\end{prop}

\begin{proof}
Since a map between $l$-adic
$\Gal(F^{\sep}/F)$-representations with different weights is zero,
we have $d_r^{i,\,j} = 0$ for $r \geq 2$.
Hence the weight spectral sequence degenerates at $E_2$.
\end{proof}

\begin{rem}
Proposition \ref{Degeneracy_weak} is a special case of
the first assertion of Theorem \ref{MainTheorem}
(see also \cite{Rapoport-Zink}, Satz 2.10,
\cite{CNakayama}, Theorem 0.1).
In \S \ref{SectionDegeneration}, we will prove the general case
without any assumption on residue fields
by a specialization argument.
\end{rem}

\subsection{A conjecture}

The following conjecture is also called
the weight-monodromy conjecture.
Sometimes, this is called {\it Deligne's conjecture on the purity
of monodromy filtration} in the literature.

\begin{conj}[\cite{HodgeI}, \cite{WeilII}, see also \cite{Rapoport-Zink},
\cite{Rapoport}, \cite{Illusie1}, \cite{Illusie2}]
\label{WMC_filtration}
Let notation be as in the beginning of this section.
Recall that
$X_K$ is a proper smooth variety of dimension $n$ over $K$
and $V := H^w_{\text{\rm \'et}}(X_{K^{\sep}},\Q_l)$.
Assume that the residue field $F$ is a purely inseparable extension of
a finitely generated extension of a prime field.
Let $M$ be the monodromy filtration on $V$
(\S \ref{SubsectionMonodromyFiltration}),
and $W$ the weight filtration on $V$
(\S \ref{SubsectionWeightFiltration}).
Then, we have $M_k V = W_{w+k} V$ for all $k$.
\end{conj}

\begin{rem}
\label{KnownCasesWMC}
Conjecture \ref{WMC_filtration} is known to hold
if either $X_K$ is an abelian variety or $\dim X_K \leq 2$
(\cite{SGA7-I}, \cite{Rapoport-Zink}, \cite{deJong}, \cite{TSaito}).
In mixed characteristic and in dimension $\geq 3$,
Conjecture \ref{WMC_filtration} is open in general.
The author proved Conjecture \ref{WMC_filtration}
for certain threefolds (\cite{Ito2}) and $p$-adically uniformized
varieties (\cite{Ito3}).
In characteristic $p>0$, if $X_K/K$ is the henselization
of a family of varieties over a curve over a finite field,
Deligne proved Conjecture \ref{WMC_filtration}
in his proof of the Weil conjectures
(\cite{WeilII}, Th\'eor\`eme 1.8.4).
Note that, in \S \ref{SubsectionApplicationEqualChar},
we will prove the equal characteristic case of
Conjecture \ref{WMC_filtration} by reducing to the case
treated by Deligne.
\end{rem}

Note that we can state Conjecture \ref{WMC_filtration}
for a proper smooth variety $X_K$ over $K$
with an assumption on the residue field $F$.
On the other hand, we can state Conjecture \ref{WMC}
for a proper strictly semistable scheme $\X$ over $\O_K$
without any assumption on $F$.
The relation between these two conjectures is as follows
(see also Remark \ref{RelationTwoConjectures} below).

\begin{prop}
\label{WMC_equivalence}
Let notation be as in the beginning of this section.
Assume that the residue field  $F$ is
a purely inseparable extension of
a finitely generated extension of a prime field,
and there exists a proper strictly semistable
scheme $\X$ over $\O_K$ such that the generic fiber
$\X_K$ of $\X$ is isomorphic to $X_K$.
Then, Conjecture \ref{WMC_filtration} for $X_K$ and
Conjecture \ref{WMC} for $\X$ are equivalent to each other.
\end{prop}

\begin{proof}
As we saw in \S \ref{SubsectionApplicationWeightSpectralSequence},
the weight spectral sequence degenerates at $E_2$
(Proposition \ref{Degeneracy_weak}), and it induces
the weight filtration $W$ on $V$ by the Weil conjectures.
The monodromy operator $N$ in \S \ref{SubsectionMonodromyFiltration}
is induced from the monodromy operator $N$ on the weight
spectral sequence (\S \ref{SectionIntroduction}).
Therefore, according to the definition of the monodromy filtration
(\S \ref{SubsectionMonodromyFiltration}),
Conjecture \ref{WMC_filtration} holds for $X_K$
if and only if $N^r \colon E_2^{-r,\,w+r}(r) \to E_2^{r,\,w-r}$
is an isomorphism for all $r,w$.
\end{proof}

\begin{rem}
\label{RelationTwoConjectures}
If we vary $K$,
Conjecture \ref{WMC} and Conjecture \ref{WMC_filtration}
are equivalent to each other.
Precisely speaking,
by using the geometric construction in \S \ref{SectionConstruction},
we can prove that Conjecture \ref{WMC_filtration}
implies Conjecture \ref{WMC} (see Remark \ref{MixedCharacteristicCase}).
Conversely, by de Jong's alteration,
we can also prove that Conjecture \ref{WMC}
implies Conjecture \ref{WMC_filtration}
(see Remark \ref{RemarkEqualChar}).
\end{rem}

\section{Some lemmas}

Here we list some lemmas which immediately
follow from the basic properties of the weight spectral sequences.
Let notation be as in \S \ref{SectionIntroduction}.
Recall that $\X$ is a proper strictly semistable scheme over
a henselian discrete valuation ring $\O_K$.
In this section, we do not put any assumption
on the residue field $F$.

\begin{lem}
\label{Lemma1}
The weight spectral sequence
$E_1^{i,\,j} \Rightarrow E^{i+j}$ of $\X/\O_K$
degenerates at $E_2$ if and only if the following equality holds:
$$ \sum_{i,j} \dim_{\Q_l} E_2^{i,\,j} = \sum_{k} \dim_{\Q_l} E^k. $$
\end{lem}

\begin{proof}
It is easy to see that the above equality holds if and only if
all maps $d_r^{i,\,j}$ for $r \geq 2$ are zero.
Hence we have the assertion.
\end{proof}

\begin{lem}
\label{Lemma2}
Conjecture \ref{WMC} depends only on the geometric special fiber of $\X$
in the following sense. Namely, if $\X/\O_{K}$ (resp.\ $\X'/\O_{K'}$)
is a proper strictly semistable scheme over
a henselian discrete valuation ring $\O_K$ (resp.\ $\O_{K'}$)
such that there exists a separably closed field $F$
containing the residue fields of $\O_K$ and $\O_{K'}$,
and $\X \otimes_{\O_K}\!F \cong \X' \otimes_{\O_{K'}}\!F$,
then Conjecture \ref{WMC} for $\X/\O_{K}$ and $\X'/\O_{K'}$
are equivalent to each other.
\end{lem}

\begin{proof}
In general, the $E_1$-terms of the weight spectral sequence
and the maps $d_1^{i,\,j},\ N$ on them are explicitly written
in terms of the geometric special fiber
(for details, see \cite{Rapoport-Zink}, Satz 2.10).
Hence, the $E_1$-terms and the maps $d_1^{i,\,j},\ N$ on them
are the same for $\X/\O_K,\ \X'/\O_{K'}$.
The maps in question
$N^r \colon E_2^{-r,\,w+r}(r) \to E_2^{r,\,w-r}$
are also the same for $\X/\O_K,\ \X'/\O_{K'}$.
Therefore, we have the assertion.
\end{proof}

\begin{lem}
\label{Lemma3}
Let $K'$ be a discrete valuation field
with ring of integers ${\O}_{K'}$ such that
${\O}_{K'}$ is henselian,
$K'$ is a field extension of $K$,
${\O}_{K'} \cap K = {\O}_K$,
and a uniformizer of $K$ is also a uniformizer of $K'$.
We put $\X' := \X \otimes_{\O_K}\!{\O}_{K'}$.
Then the assertions of Theorem \ref{MainTheorem}
for $\X/\O_K$ and $\X'/{\O}_{K'}$ are equivalent to each other.
\end{lem}

\begin{proof}
The left hand sides of the equality in Lemma \ref{Lemma1}
for $\X/\O_K,\ \X'/{\O}_{K'}$ are the same
since $\X/\O_K,\ \X'/{\O}_{K'}$
have the same geometric special fiber.
The same is true for the right hand sides
since the $l$-adic \'etale cohomology of the geometric generic
fibers of $\X/\O_K,\ \X'/{\O}_{K'}$
are the same (\cite{SGA4-III}, XVI, Corollaire 1.6).
Hence the assertion follows from Lemma \ref{Lemma1}
and Lemma \ref{Lemma2}.
\end{proof}

\begin{rem}
\label{RemarkCompletePerfect}
For any discrete valuation field $K$,
there is a field extension $K'/K$ such that
$K'$ is a complete discrete valuation field
with ring of integers ${\O}_{K'}$,
${\O}_{K'} \cap K = {\O}_K$,
a uniformizer of $K$ is also a uniformizer of $K'$,
and the residue field $F'$ of $K'$ is perfect
(see \cite{EGA3}, Chapitre 0, Proposition 10.3.1).
Therefore, by Lemma \ref{Lemma3},
it is enough to prove Theorem \ref{MainTheorem}
under the assumption that
$\O_K$ is complete and the residue field $F$ is perfect.
\end{rem}

\begin{rem}
We can apply Lemma \ref{Lemma3} for
the $\pi$-adic completion $\widehat{K}$ of $K$.
In this case, Lemma \ref{Lemma3} shows that
the assertions of Theorem \ref{MainTheorem} for $\X/\O_K$
depend only on the $\pi$-adic completion $\widehat{\X}$ of $\X$.
If we use the results of C.\ Nakayama in \cite{CNakayama},
we can say more. We note it for the reader's convenience
although we do not use it in this paper
(for details, see \cite{CNakayama}, \cite{Illusie2}).
Let $\X/\O_K$ be as in \S \ref{SectionIntroduction}.
C.\ Nakayama constructed
the weight spectral sequence only from the special fiber $\X_F$
with a natural log structure on it.
Therefore,
the assertions of Theorem \ref{MainTheorem} for $\X/\O_K$
depend only on the first infinitesimal neighborhood
$\X \otimes_{\O_K}\!\big(\O_K/(\pi^2)\big)$
since the log structure on $\X_F$
depends only on $\X \otimes_{\O_K}\!\big(\O_K/(\pi^2)\big)$
(e.g.\ \cite{Illusie2}, p.\ 309, the last paragraph).
\end{rem}

\section{A construction using N\'eron's blowing up}
\label{SectionConstruction}

In this section, we give a technical construction using N\'eron's blowing up.

Let $K$ be a discrete valuation field with ring of integers
$\O_K$ and residue field $F$.
Let $\pi$ be a uniformizer of $K$.
Let $\X$ be a proper strictly semistable scheme
of relative dimension $n$ over $\O_K$.
Then there exist a finite open covering $\{ U_i \}$
of $\X$, and, for each $i$, an \'etale morphism
$f_i \colon U_i \to \Spec \O_K[X_0,\ldots,X_n]/(X_0 \cdots X_{r_i}-\pi)$
over $\O_K$ for some $r_i$.

When $K$ is of equal characteristic $p>0$ (resp.\ equal characteristic $0$),
by N\'eron's blowing up (\cite{SGA7-I}, I, (0.5.2), \cite{Artin}, \S 4),
we can write $\O_K$ as a filtered inductive limit $\O_K = \varinjlim A_{\alpha}$,
where each $A_{\alpha}$ is a finitely generated
smooth $\F_p[\pi]_{(\pi)}$-algebra (resp.\ $\Q[\pi]_{(\pi)}$-algebra),
where $\F_p[\pi]_{(\pi)}$ (resp.\ $\Q[\pi]_{(\pi)}$) denotes
the localization of $\F_p[\pi]$ (resp.\ $\Q[\pi]$) at the maximal ideal
generated by $\pi$.
By a standard argument, for some $\alpha$,
there exist a scheme $\Y$ proper over $\Spec A_{\alpha}$
with a finite open covering $\{ V_i \}$, and an \'etale morphism
$g_i \colon V_i \to \Spec A_{\alpha}[X_0,\ldots,X_n]/(X_0 \cdots X_{r_i}-\pi)$
over $A_{\alpha}$ for each $i$ such that
the pullback of the tuple $(\Y/A_{\alpha},\,\{ V_i \},\,g_i)$
by $\Spec \O_K \to \Spec A_{\alpha}$
coincides with the original tuple $(\X/\O_K,\,\{ U_i \},\,f_i)$.

When $K$ is of mixed characteristic $(0,p)$,
the situation is slightly different but a similar construction exists.
By Lemma \ref{Lemma3}, we may assume that $\O_K$ is complete and
the residue field $F$ is perfect (see also Remark \ref{RemarkCompletePerfect}).
Then, by \cite{SGA7-I}, I, (0.5.3),
there exists a complete discrete valuation ring $R$ contained in $\O_K$
such that $\pi \in R$, the residue field of $R$ is
a purely inseparable extension of a finitely generated extension of $\F_p$,
and we can write $\O_K$
as a filtered inductive limit
  $\O_K = \varinjlim A_{\alpha}$,
where each $A_{\alpha}$ is a finitely generated smooth $R$-algebra.
Then, there exists a tuple $(\Y/A_{\alpha},\,\{ V_i \},\,g_i)$
with the same properties as above
such that the pullback of it by $\Spec \O_K \to \Spec A_{\alpha}$
coincides with the original tuple $(\X/\O_K,\,\{ U_i \},\,f_i)$.

In any case, $(\pi = 0)$ is a regular divisor of $\Spec A_{\alpha}$.
Let $s \in \Spec A_{\alpha}$ be the image of the closed point
of $\Spec \O_K$,
and $(A_{\alpha})_s$ the localization of $A_{\alpha}$ at $s$.
There exist elements $a_1,\ldots,a_r \in A_{\alpha}$ such that
$\{ \pi,a_1,\ldots,a_r \}$ is a regular system of parameters of
the local ring $(A_{\alpha})_s$. Let $B$ be the henselization of
the quotient $(A_{\alpha})_s/(a_1,\ldots,a_r)$.
Then, $B$ is a henselian discrete valuation ring
and the image of $\pi$ in $B$ is a uniformizer.
We put $\ZZ := \Y \otimes_{A_{\alpha}} B$.

In conclusion, we have the following cartesian diagram:
\begin{equation}
\label{KeyDiagram}
\begin{split}
  \xymatrix{
  \X \ar[d] \ar[r] & \Y \ar[d] & \ZZ \ar[l] \ar[d] \\
  \Spec \O_K \ar[r]& \Spec A_{\alpha} & \Spec B, \ar[l] }
\end{split}
\end{equation}
such that $\O_K,\,B$ are henselian discrete valuation rings,
and the images of the closed points of
$\Spec \O_K,\,\Spec B$ coincide with $s \in \Spec A_{\alpha}$.
Hence the geometric special fibers of $\X, \ZZ$ are the same.
The scheme $\Y$ is proper over $\Spec A_{\alpha}$ and
Zariski locally \'etale over
$\Spec A_{\alpha}[X_0,\ldots,X_n]/(X_0 \cdots X_{r}-\pi)$
for some $r$.
Hence $\ZZ$ is a proper strictly semistable scheme over $B$.
Finally, when
$K$ is of equal characteristic $p>0$ (resp.\ equal characteristic $0$),
the residue field of $B$ is finitely generated over $\F_p$ (resp.\ $\Q$).
When $K$ is of mixed characteristic $(0,p)$,
the residue field of $B$ is a purely inseparable extension of
a finitely generated extension of $\F_p$.

\section{Degeneration of the weight spectral sequences}
\label{SectionDegeneration}

In this section, we prove that the weight spectral sequences
degenerate at $E_2$ without any assumption on residue fields.
This result was already obtained by C.\ Nakayama
by using log geometry (\cite{CNakayama}, Theorem 0.1).
Here we give another proof without using log geometry.
We prove it by comparing the weight spectral sequences
of $\X/\O_K$ and $\ZZ/B$ in the diagram (\ref{KeyDiagram})
in \S \ref{SectionConstruction}.

\begin{prop}[\cite{CNakayama}, Theorem 0.1]
\label{PropDegeneracy}
Let notation be as in \S \ref{SectionIntroduction}.
Then, the weight spectral sequence of $\X/\O_K$
degenerates at $E_2$.
\end{prop}

\begin{proof}
Let us consider the diagram (\ref{KeyDiagram})
in \S \ref{SectionConstruction}.
By Proposition \ref{Degeneracy_weak},
the weight spectral sequence for $\ZZ/B$ degenerates at $E_2$.
Since $\X/\O_K$ and $\ZZ/B$ have the same geometric special fiber,
the left hand sides of the equality in Lemma \ref{Lemma1}
for $\X/\O_K,\ \ZZ/B$ are the same.
Hence, by Lemma \ref{Lemma1}, it is enough to show that
the right hand sides for $\X/\O_K,\ \ZZ/B$ are also the same.
By \cite{SGA4-III}, XVI, Corollaire 2.2,
this follows from the fact that
$\Y$ is proper and smooth outside $(\pi=0)$
and the images of the generic points of $\Spec \O_K, \Spec B$
do not lie on the divisor $(\pi=0)$.
\end{proof}

\section{Proof of Theorem \ref{MainTheorem}}

Here we shall prove Theorem \ref{MainTheorem}.
The first part of Theorem \ref{MainTheorem}
is already proved in Proposition \ref{PropDegeneracy}.

Let notation be as in \S \ref{SectionIntroduction}.
Assume that $K$ is of equal characteristic,
and consider the diagram (\ref{KeyDiagram})
in \S \ref{SectionConstruction}.
According to Lemma \ref{Lemma2},
it is enough to prove Conjecture \ref{WMC} for $\ZZ/B$.
Hence, we may replace $\X/\O_K$ by $\ZZ/B$.

Recall that,
if $K$ is of characteristic $p>0$ (resp.\ characteristic $0$),
$A_{\alpha}$ in the diagram (\ref{KeyDiagram})
is a finitely generated smooth $\F_p[\pi]_{(\pi)}$-algebra
(resp.\ $\Q[\pi]_{(\pi)}$-algebra),
and $\O_K$ is the henselization of the discrete valuation ring
$(A_{\alpha})_s / (a_1,\ldots,a_r)$.
Therefore, there exist a finitely generated smooth
$\F_p[\pi]$-algebra (resp.\ $\Z[\pi, 1/l]$-algebra)
$A$ contained in $\O_K$ and a scheme $f \colon \widetilde{\X} \to \Spec A$
such that $\O_K$ is the henselization of $A$ at $(\pi)$,
$(\pi = 0)$ is a regular divisor of $\Spec A$,
$f$ is proper and smooth outside $(\pi = 0)$,
and $\widetilde{\X} \otimes_A\!\O_K \cong \X$.

To prove Theorem \ref{MainTheorem},
it is enough to prove the following proposition.

\begin{prop}
\label{MainProp}
Let $\X/\O_K$ and $f \colon \widetilde{\X} \to \Spec A$ be as above.
Then, Conjecture \ref{WMC} holds for $\X/\O_K$.
\end{prop}

\begin{proof}
Let $M$ be the monodromy filtration on
$V := H^w_{\text{\rm \'et}}(\X_{K^{\sep}},\Q_l)$
as in \S \ref{SubsectionMonodromyFiltration}.
We can construct an $l$-adic sheaf version of $(V,M)$
as follows (for details, see \cite{WeilII}, Variante 1.7.8).
Since the action of $\Gal(K^{\sep}/K)$ on $V$
is tamely ramified (\cite{Rapoport-Zink}),
$V$ extends to a smooth $l$-adic sheaf on
$\varprojlim \Spec \O_K[\pi^{1/l^{m}}]$
with $\Z_l(1)$-action.
Therefore, the smooth $l$-adic sheaf $R^w f_{\ast} \Q_l$
on $(\Spec A) \setminus (\pi=0)$ extends to
a smooth $l$-adic sheaf on $\varprojlim \Spec A[\pi^{1/l^{m}}]$
with $\Z_l(1)$-action
by Zariski-Nagata's purity (\cite{SGA2}, X, Th\'eor\`eme 3.4).
There is a natural section of
$\varprojlim \Spec A[\pi^{1/l^{m}}] \to \Spec A$
on $(\pi=0)$ because
$\Spec A[\pi^{1/l^{m}}] \to \Spec A$ is totally ramified along $(\pi=0)$.
By pulling back, we have an $l$-adic sheaf ${\mathcal F}$
on $(\pi=0)$ with $\Z_l(1)$-action.
By the same way as in \S \ref{SubsectionMonodromyFiltration},
we have a filtration ${\mathcal M}$ on ${\mathcal F}$
such that the pair $({\mathcal F}, {\mathcal M})$
specializes to $(V, M)$ in an obvious way.

Since $A$ is a smooth $\F_p[\pi]$-algebra (resp.\ $\Z[\pi, 1/l]$-algebra),
for a closed point $t \in (\pi=0) \subset \Spec A$,
there exists a curve $C \subset \Spec A$ over a finite field
which passes through $t$ and intersects transversally with $(\pi=0)$.
Then, the filtration ${\mathcal M}$ on ${\mathcal F}$
specializes to the monodromy filtration on the restriction of
$R^w f_{\ast} \Q_l$ to $C \setminus (\pi=0)$ at $t \in C$
in the sense of \cite{WeilII}, (1.7.2).
Since $C$ is a curve over a finite field,
the stalk $(\Gr_k^{\mathcal M} {\mathcal F})_t$ has weight $w+k$,
and the smooth $l$-adic sheaf
$\Gr_k^{\mathcal M} {\mathcal F}$ on $(\pi=0)$ has weight $w+k$
(Remark \ref{KnownCasesWMC}, see also 
\cite{WeilII}, Th\'eor\`eme 1.8.4, Corollaire 1.8.7).
Since $F$ is the function field of the divisor $(\pi=0)$ and
$\Gr_k^{M} V$ is the stalk of $\Gr_k^{\mathcal M} {\mathcal F}$
at the generic point of $(\pi=0)$,
Conjecture \ref{WMC_filtration} holds for $\X_K$
(see also our definition of the notion of weights
in \S \ref{SubsectionWeightFiltration}).
Therefore, by Proposition \ref{WMC_equivalence},
Conjecture \ref{WMC} holds for $\X$.

Hence the assertion of Proposition \ref{MainProp} is proved
and the proof of Theorem \ref{MainTheorem} is complete.
\end{proof}

\begin{rem}
\label{MixedCharacteristicCase}
The above proof is valid only in the equal characteristic case.
By the same argument as above, we can prove that
Conjecture \ref{WMC_filtration} in mixed characteristic $(0,p)$
implies Conjecture \ref{WMC} in mixed characteristic $(0,p)$.
However, we can not prove Conjecture \ref{WMC} in mixed characteristic.
The problem is as follows.
In mixed characteristic $(0,p)$,
$A$ is a finitely generated smooth $R$-algebra,
where $R$ is a complete discrete valuation ring of mixed characteristic $(0,p)$
as in \S \ref{SectionConstruction}.
Since $R$ is of mixed characteristic $(0,p)$,
all characteristic $p>0$ points of $\Spec A$
are lying over the closed point of $\Spec R$.
Therefore, all curves $C \subset \Spec A$ over a finite field
are contained in the divisor $(\pi=0)$,
and there is no curve $C \subset \Spec A$ over a finite field
intersecting transversally with $(\pi=0)$.
\end{rem}

\section{Some consequences}
\label{SectionConsequence}

\subsection{The equal characteristic case of Conjecture \ref{WMC_filtration}}
\label{SubsectionApplicationEqualChar}

By a similar specialization argument,
we can also prove the equal characteristic case of
Conjecture \ref{WMC_filtration}
(for the finite residue field case, see also \cite{Terasoma}).

\begin{prop}
\label{EqualCharWMC}
Let notation be as in the beginning of \S \ref{SectionFiltrations}.
If $K$ is of equal characteristic and $F$ is a purely inseparable extension of
a finitely generated extension of a prime field,
then Conjecture \ref{WMC_filtration} holds.
\end{prop}

\begin{proof}
By Grothendieck's monodromy theorem,
after replacing $K$ by a finite extension of it,
the action of $\Gal(K^{\sep}/K)$ on
$V := H^w_{\text{\rm \'et}}(X_{K^{\sep}},\Q_l)$
is tamely ramified and the action of
the inertia group $I_K$ on $V$ is unipotent
(see \S \ref{SubsectionMonodromyFiltration}).
By the same way as in \S \ref{SectionConstruction},
if $K$ is of characteristic $p>0$ (resp.\ characteristic 0),
there exist a finitely generated smooth $\F_p[\pi]$-algebra
(resp.\ $\Z[\pi, 1/l]$-algebra) $A$ contained in $\O_K$
and a scheme $f \colon \widetilde{\X} \to \Spec A$
such that $f$ is proper and smooth outside $(\pi=0)$
and $\widetilde{\X} \otimes_A K \cong X_K$.
Here we do not make any assumption on the fiber of $f$ over $(\pi=0)$.
Then, by the same way as in the proof of Proposition \ref{MainProp},
there exist a smooth $l$-adic sheaf ${\mathcal F}$ on $(\pi=0)$
and a filtration ${\mathcal M}$ on it
such that the pullback of $({\mathcal F}, {\mathcal M})$
by $\Spec \O_K \to \Spec A$ is $(V,M)$,
and $\Gr_k^{\mathcal M}{\mathcal F}$ has weight $w+k$
in the usual sense.
Since $F$ is an extension of the residue field $\kappa(s)$
at some point $s \in (\pi=0)$, $\Gr_k^M V$ has weight $w+k$
in the sense of \S \ref{SubsectionWeightFiltration}.
Hence Conjecture \ref{WMC_filtration} holds for $X_K$.
\end{proof}

\begin{rem}
\label{RemarkEqualChar}
It is also possible to deduce Proposition \ref{EqualCharWMC}
from Theorem \ref{MainTheorem} as follows.
By de Jong's alteration (\cite{deJong}, Theorem 6.5),
after replacing $K$ by a finite extension of it,
there is a proper strictly semistable scheme $\Y$ over $\O_K$,
and a proper surjective generically finite morphism
$g \colon \Y_K \to X_K$.
We put $V := H^w_{\text{\rm \'et}}(X_{K^{\sep}},\Q_l)$
and $V' := H^w_{\text{\rm \'et}}(\Y_{K^{\sep}},\Q_l)$.
Since $\Y_K, X_K$ are proper and smooth over $K$,
there is a $\Gal(K^{\sep}/K)$-equivariant trace map
$g_{\ast} \colon V' \to V$.
We know that the composite $g_{\ast} \circ g^{\ast}$ is
the multiplication-by-$d$ map on $V$,
where $d$ denotes the separable degree of $g$
at the generic point of $X_K$.
Hence $g^{\ast} \colon V \to V'$ is injective and
$V$ is a direct summand of $V'$ as $\Gal(K^{\sep}/K)$-representations.
Therefore, the restriction of the monodromy (resp.\ weight) filtration
on $V'$ to $V$ is the monodromy (resp.\ weight) filtration on $V$.
By Theorem \ref{MainTheorem} and Proposition \ref{WMC_equivalence},
Conjecture \ref{WMC_filtration} holds for $\Y_K$.
Hence Conjecture \ref{WMC_filtration} holds also for $X_K$.
\end{rem}

\subsection{A Hodge analogue over $\C$}
\label{SubsectionApplicationHodge}

When $F=\C$, Steenbrink originally considered
a Hodge analogue of Conjecture \ref{WMC} (\cite{Steenbrink}).
As an application of Theorem \ref{MainTheorem},
we can give a modulo $p>0$ reduction proof of it.

\begin{prop}[\cite{Steenbrink}, Proposition 5.14]
\label{HodgeAnalogue}
Let notation be as in \cite{Steenbrink}, \S 5.
Then, Schmid's limit Hodge structures coincide with
Steenbrink's limit Hodge structures.
\end{prop}

\begin{proof}
It is enough to prove \cite{Steenbrink}, Theorem 5.9
which is nothing but a Hodge analogue of Conjecture \ref{WMC}.
Hence, by the comparison theorem between \'etale
and singular cohomology (\cite{SGA4-III}, XI, Th\'eor\`eme 4.4),
this follows from Theorem \ref{MainTheorem}.
\end{proof}

\begin{rem}
Note that, Steenbrink's original argument in \cite{Steenbrink}, Theorem 5.9
is incomplete (see \cite{MSaito1}, 4.2.5, \cite{SaitoZucker}, \S 2.3,
\cite{GuillenNavarroAznar}).
\end{rem}

\begin{rem}
Conversely, by Lefschetz principle, if $K$ is of equal characteristic $0$,
Conjecture \ref{WMC} follows from the results over $\C$.
Therefore, we have two completely different proofs of Conjecture \ref{WMC}
in the equal characteristic $0$ case.
Illusie pointed out to the author that there is a technical difference
between them. In the Hodge theoretic proof, the proper case is reduced to
the projective case by Chow's lemma because polarized Hodge structures
are used (see \cite{SaitoZucker}, Remarks (0.5), (iv)).
On the other hand, in this paper,
we can directly treat the proper case because
Deligne's results in \cite{WeilII} are valid
without the projectivity assumption.
\end{rem}


\begin{thebibliography}{GNA}

\bibitem[Ar]{Artin} M. Artin, {\it Algebraic approximation of
  structures over complete local rings},
  Inst. Hautes \'Etudes Sci. Publ. Math. No. 36, (1969), 23--58.
\bibitem[De1]{HodgeI} P. Deligne, {\it Th\'eorie de Hodge I},
  in {\it Actes du Congr\`es International des Math\'ematiciens
  (Nice, 1970), Tome 1}, 425--430, Gauthier-Villars, Paris, 1971.
\bibitem[De2]{WeilI} P. Deligne, {\it La conjecture de Weil I},
  Inst. Hautes \'Etudes Sci. Publ. Math. No. 43, (1974), 273--307.
\bibitem[De3]{WeilII} P. Deligne, {\it La conjecture de Weil II},
  Inst. Hautes \'Etudes Sci. Publ. Math. No. 52, (1980), 137--252.
\bibitem[GH]{GoertzHaines} U. G\"ortz, T. Haines,
  {\it The Jordan-H\"older series for nearby cycles on some
  Shimura varieties and affine flag varieties},
  preprint, math.AG/0402143, 2004.
\bibitem[GNA]{GuillenNavarroAznar}
  F. Guill\'en, V. Navarro Aznar,
  {\it Sur le th\'eor\`eme local des cycles invariants}
  Duke Math. J. {\bf 61} (1990), no. 1, 133--155.
\bibitem[Il1]{Illusie1} L. Illusie, {\it Autour du th\'eor\`eme de
  monodromie locale}, P\'eriodes $p$-adiques (Bures-sur-Yvette, 1988),
  Ast\'erisque No. 223, (1994), 9--57.
\bibitem[Il2]{Illusie2} L. Illusie,
  {\it An overview of the work of K.\ Fujiwara, K.\ Kato, and
  C.\ Nakayama on logarithmic etale cohomology},
  in {\it Cohomologies $p$-adiques et applications arithm\'etiques, II},
  Ast\'erisque No. 279 (2002), 271--322.
\bibitem[It1]{Ito1} T. Ito, {\it Weight-monodromy conjecture over
  positive characteristic local fields}, Master's thesis,
  University of Tokyo, 2001.
\bibitem[It2]{Ito2} T. Ito, {\it Weight-monodromy conjecture for
  certain threefolds in mixed characteristic},
  Internat. Math. Res. Notices {\bf 2004}, no.~2, 69--87.
\bibitem[It3]{Ito3} T. Ito, {\it Weight-monodromy conjecture for
  $p$-adically uniformized varieties}, math.NT/0301201, 2003,
  to appear in Invent. Math.
\bibitem[dJ]{deJong} A. J. de Jong, {\it Smoothness, semi-stability
  and alterations}, Inst. Hautes \'Etudes Sci. Publ. Math. No. 83,
  (1996), 51--93.
\bibitem[Na]{CNakayama} C. Nakayama, {\it Degeneration of $l$-adic
  weight spectral sequences},
  Amer. J. Math. {\bf 122} (2000), no.~4, 721--733.
\bibitem[Ra]{Rapoport} M. Rapoport,
  {\it On the bad reduction of Shimura varieties},
  in {\it Automorphic forms, Shimura varieties, and $L$-functions,
  Vol.\ II (Ann Arbor, MI, 1988)}, 253--321, Academic Press, Boston,
  MA, 1990.
\bibitem[RZ]{Rapoport-Zink} M. Rapoport, T. Zink, {\it \"{U}ber die
  lokale Zetafunktion von Shimuravariet\"aten. Monodromiefiltration
  und verschwindende Zyklen in ungleicher Charakteristik},
  Invent. Math. {\bf 68} (1982), no.~1, 21--101.
\bibitem[SaM]{MSaito1} M. Saito, {\it Modules de Hodge polarisables},
  Publ. Res. Inst. Math. Sci. {\bf 24} (1988), no.~6, 849--995 (1989).
\bibitem[SaZ]{SaitoZucker} M. Saito, S. Zucker,
  {\it The kernel spectral sequence of vanishing cycles},
  Duke Math. J. {\bf 61} (1990), no. 2, 329--339.
\bibitem[SaT]{TSaito} T. Saito, {\it Weight spectral sequences and
  independence of $\ell$}, 
  J. Inst. Math. Jussieu {\bf 2} (2003), 583--634.
\bibitem[Sc]{Schmid} W. Schmid, {\it Variation of Hodge structure :
  the singularities of the period mapping},
  Invent. Math. {\bf 22} (1973), 211--319.
\bibitem[Se]{Serre} J.-P. Serre, {\it Corps locaux}, Deuxieme edition,
  Hermann, Paris, 1968.
\bibitem[ST]{SerreTate} J.-P. Serre, J. Tate, {\it Good reduction of
  abelian varieties}, Ann. of Math. (2) {\bf 88} (1968), 492--517.
\bibitem[St]{Steenbrink} J. Steenbrink, {\it Limits of Hodge structures},
  Invent. Math. {\bf 31} (1975/76), no.~3, 229--257.
\bibitem[Te]{Terasoma} T. Terasoma, {\it Monodromy weight filtration
  is independent of $l$}, preprint, math.AG/9802051, 1998.
\bibitem[EGA3]{EGA3} {\it \'El\'ements de g\'eom\'etrie alg\'ebrique. III.
  \'Etude cohomologique des faisceaux coh\'erents. I.},
  Inst. Hautes \'Etudes Sci. Publ. Math. No. 11 (1961), 167 pp.
\bibitem[SGA2]{SGA2} {\it Cohomologie locale des faisceaux coh\'erents et
  th\'eor\`emes de Lefschetz locaux et globaux}, Adv. Stud. Pure Math., 2,
  North-Holland, Amsterdam, 1968.
\bibitem[SGA4-III]{SGA4-III} {\it Th\'eorie des topos et cohomologie \'etale
  des sch\'emas. Tome 3}, Lecture Notes in Math., 305, Springer, Berlin, 1973.
\bibitem[SGA7-I]{SGA7-I} {\it Groupes de monodromie en g\'eom\'etrie
  alg\'ebrique. I}, Lecture Notes in Math., 288, Springer, Berlin, 1972.
\end{thebibliography}
\end{document}